\documentclass[12pt]{amsart}
\usepackage{amsmath}
\usepackage{amssymb}
\newtheorem{theorem}{Theorem}
\newtheorem{remark}{Remark}
\newtheorem{corollary}{Corollary}
\newtheorem{proposition}[theorem]{Proposition}
\newtheorem{lemma}[theorem]{Lemma}

\newenvironment{proof*}{\vskip 2mm\noindent {}}{\hfill $\Box$ \vskip 2mm}

\newcommand{\B}{\mathbb B}
\newcommand{\C}{\mathbb C}
\newcommand{\D}{\mathbb D}
\newcommand{\R}{\mathbb R}
\newcommand{\eps}{\varepsilon}

\def\grad{\operatorname{grad}}
\def\Re{\operatorname{Re}}
\def\Im{\operatorname{Im}}

\def\span{\operatorname{span}}
\long\def\comment#1{}
\begin{document}

\title{Two-dimensional slices of non-pseudoconvex open sets}

\author{Nikolai Nikolov}
\address{Institute of Mathematics and Informatics\\ Bulgarian Academy of
Sciences\\1113 Sofia, Bulgaria} \email{nik@math.bas.bg }

\author{Peter Pflug}
\address{Carl von Ossietzky Universit\"at Oldenburg\\
Institut f\"ur Mathematik, Fakult\"at V\\ Postfach 2503\\ D-26111
Oldenburg, Germany}\email{peter.pflug@uni-oldenburg.de}



\begin{thanks}{This paper was written during the stay of the first-named author
at the Carl von Ossietzky Universit\"at, Oldenburg (February--March 2010) supported by a
DFG grant 436Pol113/106/2. The authors would like to thank Pascal J.~Thomas for helpful discussions
about Proposition \ref{sli} and for pointing out a serious mistake in the proof of Lemma \ref{3} in a former version of this paper.}
\end{thanks}

\begin{abstract} Let $D$ be a non-pseudoconvex open set in $\C^3$ and $S$ be the union of
all two-dimensional planes with non-empty and non-pseudoconvex intersection with $D.$ Sufficient
conditions are given for $\C^3\setminus S$ to belong to a complex line. Moreover, in the $\mathcal C^2$-smooth case, it is
shown that $S=\C^n$.

\end{abstract}

\maketitle

It is well known (see \cite{Hit}, \cite{Jac}) that a open set $D\subset\C^n$ is pseudoconvex if and only if $D\cap P$ is pseudoconvex (or empty)
for every two-dimensional plane $P\subset\C^n$. Hence, pseudoconvexity is, in fact, a two-dimensional phenomenon. Let us denote by $S=S(D)$ the set of all points $a\in\C^n$ such that
there exists such a two-dimensional plane $P$, $a\in P$, such that $P\cap D$ is non-pseudoconvex. Points which do not belong to $S$ are called \textit{exceptional points for $D$}.
In \cite{Nik-Tho} the following results are shown:

\begin{enumerate}
\item  if $D\subset\C^n$, $n\geq 3$, is not pseudoconvex, then all exceptional points belong to one complex hyperplane;

\item  if $D$ is $\mathcal C^2$-smooth and non-pseudoconvex near a boundary point, then all exceptional points belong to a three-codimensional plane;

\item there is an unbounded domain in $\C^3$ with real-analytic boundary except one point, the origin, which has exactly one exceptional point, namely the origin;

\item if $D=D_0\setminus l\subset\C^3$, where $D_0$ is a pseudoconvex open set and $l$  a complex line intersecting $D_0$, then $S=\C^3\setminus(l\setminus D)$;

\item if $D$ is as before and if $l_1,l_2$ are two different complex lines with $l_j\cap D\neq\varnothing$, then $\C^3\setminus S=(l_1\cap l_2)\setminus D$.
\end{enumerate}

Looking more carefully at the example (3) which is given as
$$D=\{z\in\C^3:|z_3|^2<|z_1|^2+|z_2|^2<4|z_3|^2\},
$$
we see that the restriction of the Hessian (resp.~the Levi form) of $r(z)=|z_3|^2-|z_1|^2-|z_2|^2$
to the complex tangent hyperplane $T_p$ at any non-zero $p$ with $r(p)=0$ has eigenvalues $0,0,-1,-1$ (resp.~$0$ and $-1$).
Moreover, if  $l_p$ denotes the complex line through $0$ and $p,$ then
$T_p\setminus l_p\subset D$ and $l_p\subset\partial D.$

So, it is natural to study the Hessian as in the $\C$-convex case because here we play with intersections with complex lines.

The following proposition will play a key role in the next considerations.

\begin{proposition}\label{sli} Assume that $D$ is an open set in $\C^3$
such that $0\in\partial D$ and
$D\supset G\cap\D^3_\eps$ for some $\eps>0$ {\rm (}$\D_\eps:=\{\lambda\in\C:|\lambda|<\eps\}${\rm )}, where\footnote{For an arbitrary non-pseudoconvex
set $D$ in $\C^3,$ after an affine change of coordinates, we may assume the same but with an extra summand $c(\Im z_2)^2,$
$c<2$ (cf.~\cite{Nik-Tho}). So, after a quadratic polynomial automorphism of $\C^3,$
we have the original domain $G$ locally inscribed in $D.$}
$$
G=\{z\in\C^3:0>2\Re z_3+(\Im z_3)^2+|z_1|^2-|z_2|^2\}.
$$

If $S\neq\C^3,$ then here exists a complex line $l\in\C^2\times\{0\}\subset\C^3$ through $0$ such that
$\C^2\times\{0\}\setminus l\subset D$ near $0$, $l\subset\partial D$ near $0$, and $S\supset\C^3\setminus l.$
\end{proposition}

\begin{proof} Let $a\in\C^3.$ If $a_3\neq 0,$ then the intersection of $D$ with the complex hyperplane through the point $0,$ $a$, and
$(0,1,0)$ is non-pseudoconvex; use that Levi form of $G$ at $0$ is negative and apply the Kontinuit\"atsatz.
Hence $S\supset\{z\in\C^3:z_3\neq 0\}.$

Fix a $S\not\ni a=(a_1,a_2,0)\neq 0$ and put $l=\C a$.  It will be enough to show that
\begin{equation}\label{key}
D'\times\{0\}:=D\cap\{z\in\C^3:z_3=0\}\supset\Big(\D^2_{\eps/2}\setminus l \Big)\times\{0\}.
\end{equation}

Indeed, suppose that there is a $b=(b_1,b_2,0)\not\in S\cup l.$ It follows by (\ref{key}) that
$D'\supset\D^2_{\eps/2}\setminus\{z\in\C^2: b_1z_2=b_2z_1\}.$
Hence $0\not\in D'\supset\D^2_{\eps/2}\setminus\{0\}.$ In particular,
$D'$ is non-pseudoconvex which contradicts the fact that $a\notin S$.

To prove (\ref{key}), let $(0,0)\neq (p,q)\in\C^2$, $P=P(p,q):=\{z\in\C^3:z_3=\overline{p}(z_1-a_1)+\overline{q}(z_2-a_2)\}$, and $P_\eps=P\cap \D^3_\eps.$ Then
$$D\cap P\supset G\cap P_\eps\supset\{z\in P_\eps:|z_1+p|^2<|z_2-q|^2+s\},$$
where $$s=2\Re(\overline{p}a_1+\overline{q}a_2)+|p|^2(1-2(|a_1|+\eps)^2)-|q|^2(1+2(|a_2|+\eps)^2).$$
Thus we may find a positive $\eps'<\eps$ such that if $(p,q)\in\D^2_{\eps'},$ then
$$D\supset D_\eps:=\{z\in P:|z_1+p|^2<|z_2-q|^2+s',(z_1,z_2)\in\D^2_\eps\},$$
where
$$s'=2\Re(\overline{p}a_1+\overline{q}a_2)-|p|^{3/2}-|q|^{3/2}.$$

Take now a point $z^0=(z_1^0,z_2^0,0)\in\D^3_{\eps/2}$ with $a_1z_2^0\neq a_2z_1^0.$
Assume that $z_2^0\neq a_2$ (the case $z_1^0\neq a_1$ is similar). It is easy
to check that there is $\eta<\eps'/2$ such that for any $p\neq 0$, $|p|<\eta$, with $\arg p=\arg\left(\frac{a_1z_2^0-a_2z_1^0}{z_2^0-a_2}\right)$
and $\overline{q}=\overline{p}\frac{a_1-z_1^0}{z_2^0-a_2}$ we have $s'>0$ and $(z_1^0+p,z_2^0-q)\in\D_{\eps/2}^2$. Fix such a pair $(p,q)$ and let $P=P(p,q)$ be the corresponding plane. Observe that $z^0\in P$. Then
$$
D_\eps\supset\widetilde D_\eps:=\{z\in P:|z_1+p|^2<|z_2-q|^2+s',(z_1+p,z_2-q)\in\D^2_{\eps/2}\}.
$$
Then the envelope of holomorphy of $\widetilde D_\eps$ is
$$
H=\{z\in P:(z_1+p,z_2-q)\in\D^2_{\eps/2}\},
$$
as it is known from the study of Reinhardt domains.
Since $D\cap P$ is pseudoconvex (because $a\not\in S$),
it follows that $z^0\in H\subset D\cap P$ which proves (\ref{key}).

It remains to show that $l\subset D$ near $0.$ We know that $\hat D\times\{0\}:=D\cap\{z\in\C^3: z_3=0\}\cap\D^3_{\eps/2}$ is a pseudoconvex open set and
$\D^2_{\eps/2}\setminus l\subset \hat D\subset(\D^2_\eps)_\ast.$ By the claim below, $\hat D=\D^2_{\eps/2}\setminus l.$
Therefore it follows that $(l\cap\D^2_{\eps/2})\times\{0\}\subset\partial D.$
\smallskip

\noindent{\it Claim.} If $l=\C a$, $\|a\|=1$, is a complex line in $\C^2$ and $G\subsetneq\D^2$ is a pseudoconvex open set containing
$\D^2\setminus l,$ then $G=\D^2\setminus l.$
\smallskip

Indeed, first note that $G$ is a domain, since if $\D^2\ni\lambda a$, then for any sufficiently small $t>0$ we have that $\lambda a+t(-a_2,a_1)\in\D^2\setminus l$, where $a=(a_1,a_2)$. Obviously, we may assume that $|a_2|\leq|a_1|$ (otherwise  interchange the coordinates). Now suppose that $G\neq \D^2\setminus l$. Then there exists a point $\lambda_0 a\in G$ and so we have that \footnote{$\D_\eps(c)$ will always denote the disc with center $c\in\C$ and radius $\eps$.}
$$
V_\eps:=\D_\eps(\lambda_0a_1)\times\D_\eps(\lambda_0a_2)\subset G
$$
for a small $\eps>0$. Moreover,
simple calculation gives that
$$
(\D_\eps(\lambda_0a_1)\times\D)\cap l\subset V_\eps.
$$
Therefore, $U_\eps:=\D_\eps(\lambda_0a_1)\times\D\subset G$. Fix an arbitrary $r<1$ such that $\D_\eps(\lambda_0a_1)\subset\D_r$. Then
$$
F(z):=\frac{1}{2\pi i}\int_{|\eta|=(1+r)/2}\frac{f(z_1,\eta)}{z_2-\eta}d\eta
$$
defines a holomorphic function on $W_r:=\D_r\times\D_{(1+r)/2}$, which coincides with $f$ on $\D_\eps(\lambda_0a_1)\times\D_{(1+r)/2}$. Noting that $W_r\cap G$ is connected, we get $f=F$ on $W_r\cap G$. With $r\nearrow 1$ we finally get a holomorphic extension of $f$ to the whole bidisc; a contradiction.
\end{proof}

We like to mention that the claim before remains true for any complex (not necessarily passing through the origin) line intersecting the bidisc (with a similar proof).

Next we state the first main result.

\begin{theorem}\label{2} Let $D$ be an open set in $\C^3$ and $p$ one of its boundary points. Assume that
$D$ has a $\mathcal C^1$-defining function $r$ near $p$ (i.e.~its gradient at $p$ is non-zero) which is twice differentiable at $p$
and that the Hessian of $r$ at $p$ restricted on the complex tangent hyperplane $T_p$ at $p$
has two negative eigenvalues and one non-positive eigenvalue.

If $S\neq\C^3,$ then there exists a complex line $l$ through $p$ such that $T_p\setminus l\subset D$ near $p,$
$l\subset\partial D$ near $p$ and $S\supset\C^3\setminus l.$

Moreover, for any point $b\not\in l$ or $b=p,$ there exists a
two-dimensional complex plane $P_b\subset\C^3$ through $p$ and $b$
such that the intersection $D\cap P_b$ is non-pseudoconvex (so $p\in S$).
On the other hand, this fails to be true if $b\in l$ and $b\neq p.$

In particular, the Hessian (and hence the Levi form) at $p$ vanishes along $l;$ therefore, the Hessian (resp.~the Levi form) at $p$
restricted to $T_p$ has two (resp.~one) negative and two (resp.~one) zero eigenvalues.
\end{theorem}

Before presenting the proof we like to mention that the used conditions for the eigenvalues are optimal. Indeed, there are pseudoconvex examples with one negative and two zero eigenvalues, as well as
with two negative eigenvalues; for example,
$$\{z\in\C^3:\Re z_1+2(\Re z_2)^2-(\Im z_2)^2<0\} \quad\text{ and }
$$
$$
\{z\in\C^3:\Re z_1+2(\Re z_2)^2-(\Im z_2)^2+2(\Re z_3)^2-(\Im z_3)^2<0\}.
$$

\begin{proof} We may assume that $p=0$ and $T_0=\{z\in\C^3: z_3=0\}.$ Take a point $a\in\C^3.$ If $a_3\neq 0,$
then the beginning of the proof of Proposition \ref{sli} provides a respective complex hyperplane through $a$ (and even through $0$) (because the Levi form at $0$ has a negative eigenvalue; see the arguments below).  Let $a_3=0.$ Let $v_1,v_2,v_3,v_4$ be the (orthonormal) eigenvectors of respective matrix to the restriction $H$ of the Hessian at $0$ on $z_3=0$ which correspond to the eigenvalues $\lambda_1<0,\lambda_2<0,\lambda_3\le 0,\lambda_4.$ Consider the complex line $0\in s\subset L:=\span(v_1,v_2,v_3)$ (if $v_4=(t_1,t_2)\in\C^2,$ then $s=\{z\in\C^2:{\overline t_1}z_1=-{\overline t_2}z_2\}$). We may assume that $s$ is the $z_1$-line.
If $H$ is negative on $s$, then we may use Proposition \ref{sli} (we have the respective inscribed model $G$). Otherwise, $\lambda_3=0$ and, after a rotation in $s,$ we may assume that $\span(v_3)=\Im z_1,$ i.e.~$H=H(z_1,z_2)$ does not depend on $\Im z_1.$ We shall show that $S=\C^3.$ After a rotation in the $z_2$-line, we may assume that $L$ contains the
$\Im z_2$-axis. Then $H(\Re z_1,\Im z_1,0,\Im z_2)\le c_1(|\Re z_1|^2+|\Im z_2|^2),$ where $c_1=\max(\lambda_1,\lambda_2).$ Fix a $c_2>c_1.$ Then we may find a $c_3\gg 1$ such that
\begin{multline*}
D_3\times\{0\}:=D\cap\{z\in\C^3:z_3=0\}\supset\\ \Big(\{z\in\C^2:c_3(\Re z_2)^2+c_2(\Im z_2)^2<-c_2(\Re z_1)^2+o(|\Im z_1|^2)\}\Big)\times\{0\}.
\end{multline*}
So, for any $\eps>0$ there is a bidisc $\D^2_\delta$ such that
$$D_3\supset\{z\in\D^2_\delta:c_3(\Re z_2)^2+c_2(\Im z_2)^2<-c_2(\Re z_1)^2-\eps(\Im z_1)^2\}.$$
After a dilation of coordinates, we may find $c>0$ such that for any $\eps>0$ there is a bidisc
$\D^2_\delta$ with
$$D_3\supset\{z\in\D^2_\delta:(\Re z_2)^2-c(\Im z_2)^2<(\Re z_1)^2-\eps(\Im z_1)^2\}.$$

The following lemma shows that $D_3$ is non-pseudoconvex which completes the proof.

\begin{lemma}\label{3} Let
$$
G:=\{z=(x_1+iy_1,x_2+iy_2)\in\D^2: x_2^2-cy_2^2<x_1^2-\eps y_1^2\}
$$
with $c>\max(\eps,0)$. Then any open set $D\subset\C^2$ satisfying $G\subset D$ and $0\in\partial D$ is non-pseudoconvex.
\end{lemma}

\begin{proof} Note that, decreasing $c$ and increasing $\eps$, we shrink $G.$ So, we may assume that
either $0<\eps<c<1$ or $1<\eps<c.$

Let first $0<\eps<c<1.$ Choose $1<p<q$ such that $\eps=\frac{p-1}{p+1}$ and
$c=\frac{q-1}{q+1}.$ Moreover, let $\Phi:\C^2\rightarrow\C^2$ be the following biholomorphic mapping: $\Phi(z):=(\sqrt{\frac{p}{p+1}}z_1,\sqrt{\frac{q}{q+1}}z_2)$.
Then a simple calculation shows that
$$
\Phi(D)\supset\Phi(G)=\{w\in\D_{\sqrt{\frac{p+1}{p}}}\times\D_{\sqrt{\frac{q+1}{q}}}:\Re(w_2^2-w_1^2)<|w_1|^2/p-|w_2|^2/q\}.
$$
So, we may suppose that
$$\Phi(D)\supset \tilde G=\{w\in\D_{2\eta}^2:\Re(w_2^2-w_1^2)<|w_1|^2/p-|w_2|^2/q\},$$
where $\eta$ is suitably chosen.

Put $S_t=\{w\in\D^2_\eta: w_2^2-w_1^2=t\}$ and $T_t=\{w\in\partial D_\eta:w_2^2-w_1^2=t\}$, where $t<0$ ($S_t\neq\varnothing$ for $t\sim 0$). Then the triangle inequality shows that
$S_t\subset\tilde G\subset\Phi(D)$ and $\bigcup_{0<-t\ll\eta^2} T_t\Subset\tilde G\subset\Phi(D).$ Moreover, note that $S_t\ni (\sqrt{-t},0)\underset{t\nearrow 0}\rightarrow (0,0)\in\partial\Phi(D)$. Recall that by the maximum principle on analytic hypersurfaces we have $\|f\|_{S_t\cup T_t}=\|f\|_{T_t}$ for any function $f$ holomorphic on $\Phi(D)$ ($-t\ll 1$). Applying the Kontinuit\" atssatz implies that $\Phi(D)$ cannot be pseudoconvex. Therefore, $D$ is non-pseudoconvex.

It remains to note that the case $1<\eps<c$ reduces to the case considered before by the biholomorhic mapping
$\Psi(z)=(i\sqrt{\eps}z_2,i\sqrt{c}z_1)$ (and replacing  $\D^2$ by $\Psi(\D^2)$).
\end{proof}
\end{proof}

\begin{remark}{\rm The assumption $c>\max(\eps,0)$ in Lemma \ref{3} is essential; for example, if $\eps\ge c=1,$ then $G$ is pseudoconvex,
and if $c\le 0,$ then $G\subset\D_\ast\times\D$ and the last domain is pseudoconvex.}
\end{remark}

Note that Proposition 1 also implies the following

\begin{proposition}\label{4}  Let $D$ and $r$  be as in Theorem \ref{2} and assume that the restriction of the Hessian of $r$ at
$p$ on a complex line in $T_p$ has two negative eigenvalues. Then the same conclusions hold as in Theorem \ref{2}.
\end{proposition}

Indeed, we may assume that $p=0$ and the Hessian is negative on the $z_2$-line. Then, after a dilation of coordinates,
we have the inscribed model $G$ as in Proposition \ref{sli}.
\smallskip

Finally, we have the following result.

\begin{proposition}\label{5}  Let $D,$ $r$ and $T_p$ be as in Theorem \ref{2}. Assume that the Levi form of $r$ at $p$ restricted on $T_p$ has two negative eigenvalues. Then $S=\C^3.$ Moreover, for any point $a\in C^3$ there exists a
two-dimensional plane $P\subset\C^3$ containing $p$ and $a,$ such that the intersection
$D\cap P$ is non-pseudoconvex for every two-dimensional plane $P\subset\C^n$,
\end{proposition}

\begin{proof} We may assume that $p=0$ and $T_0=\{z\in\C^3:z_3=0\}.$ If $a\not\in T_0,$ then $a\in S$ (as before).
To see that $a\in S$ if $a\in T_0,$ it suffices to show that $D'=D\cap T_0$ (considered as two-dimensional set) is non-pseudoconvex. Note that
$$D'\cap\D_\eta^2=\{z\in\D_\eta^2: \tilde r(z)=r(z,0)=c\Re q(z)+\mathcal L(z)+o(|z|^2)<0\},$$ where $\mathcal L$ denotes the  Levi form of $\tilde r$ at $0$,
$c\in\{0,1\}$, and $q$ is a (holomorphic) homogenous polynomial of order 2. The case $c=0$ follows immediately by Hartog's theorem. So let $c=1$. Then there exist $\eps\in(0,\eta)$ and $d>0$ such that
$r(z)\leq \Re q(z)-d\|z\|^2$, $z\in\D_{\eps}^2$.
Put $A_t=\{z\in\D^2_\eps:q(z)=t\}.$ Then there exist a sequence $(a_j)$ with $a_j\rightarrow 0$ such that $t_j:=q(a_j)=\Re q(a_j)<0$ and $A_{t_j}$ is an analytic hypersurface (use that $q$ is an open mapping). Put $S_j:=A_{t_j}$ and $T_j:=A_{t_j}\cap\partial\D^2_\eps.$ Then one may find a $k$ such that $\bigcup_{j=k}^\infty S_j\subset D'$ and
$\bigcup_{j=k}^\infty T_j\Subset D'$ (use the upper estimate for $r$ from above).
Suppose now that $D'$ is pseudoconvex.  Then the Kontinuit\"atsatz implies (similar as in the proof of Lemma \ref{3}) that $\bigcup_{j=k}^\infty S_j\Subset D'$.
On the other hand, we know that $S_j\ni a_j\to 0\in\partial D'$; a contradiction. Hence $D'$ is non-pseudoconvex.
\end{proof}

In contrast to (3) we have the following general result.

\begin{theorem}\label{add} Let $D\subset\C^n$ be a $\mathcal C^1$-smooth open set having a $\mathcal C^2$-smooth
non-pseudoconvex boundary point. If the set of non $\mathcal C^2$-smooth boundary points of $D$ has zero
$(2n-2)$-Hausdorff measure, then $D$ has no exceptional points.
\end{theorem}

In particular, we may formulate the following pseudoconvex test for $\mathcal C^2$-smooth open sets.

\begin{corollary}  A $\mathcal C^2$-smooth open set $D$ in $\C^n$ is pseudoconvex if and only if it allows one exceptional point.
\end{corollary}

\begin{proof}[Proof of Theorem \ref{add}] Put $\partial _gD$ the set of all $\mathcal C^2$-smooth boundary points of $D$. Then we find an open neighborhood $U$ of $\partial_g D$ and a defining $\mathcal C^2$-function $r$ for $\partial D$ along $\partial _g D$ on $U$ with $\|\grad r\|=1$ on $\partial_g D$.

Suppose now that there exists a point $a\in\C^n$ such that for every two-dimen\-sional plane $P$ passing through $a$ the intersection  $P\cap D$ is either empty or pseudoconvex, i.e.~$a$ is an exceptional point for $D$. We may assume that $a=0$ (use translation).

Now, take an arbitrary point $b\in\partial_g D$, $b\neq 0$, such that there exists at least one complex tangent vector $X\in T_b^\C(\partial D)$ (i.e.~$\sum_{j=1}^n\frac{\partial r}{\partial z_j}(b)X_j=0$)  such that $\mathcal Lr(b;X)<0$.

Denote by $s(z)=s_r(z)$, $z\in\partial_g D$, the smallest eigenvalue of the Levi form of $r$, restricted to the complex tangent hyperplane $T_z^\C(\partial D)$, and put
$$
t_z:=\inf\{t\in(0,1]: \forall_{\tau\in[t,1]}:\tau z\in\partial_g D \text{  and } s(\tau z)=s(z)\tau<0\}.
$$
If $t_b>0$, then we claim that $t_bb\in\partial D\setminus\partial_gD$. Indeed, we may assume that $b=(b_1,0,\dots,0)$ with $b_1=\Re b_1<0$ (use rotation).
Suppose now the contrary, namely that  $t_bb\in\partial_gD$.

So we have
$0\neq c=:t_bb\in\partial_g D$ and $s(c)=s(b)/t_0$. Since $r(tc)=0$ for $t$ belonging to some interval
$[1,1+\eps)$ we get $\partial r/\partial x_1(tc)=0$ for these $t$'s; in particular, $\partial r/\partial x_1(c)=0$.

From now we identify a complex point $z=(x_1+iy_1,\dots, x_n+iy_n)$ with its real counterpart $(x_1,\cdots,x_{2n})$ (i.e.~$y_j=x_{n+j}$ for $j=1,\dots,n$). Then
we may solve the equation $r=0$ near $c$ w.r.t., say the $k$-th coordinate ($k\neq 1$). So we have an open $\delta$-cube $Q_\delta$ with center $(c_1,0,\cdots,0)\in\R^{2n-1}$ and an interval $I:=(-\eta,\eta)$ such that
$x\in (Q\times I)\cap\partial D$ if and only if $x_k=f(\tilde x)$, $\tilde x\in Q$. We know that $f(\xi_1,\tilde 0)=0$ as long $\xi_1\in (c_1-\delta,c_1]$.

By assumption (recall that $s(c)<0$), there is  $\delta'<\delta$ a open $\delta'$-cube $V=Q_{\delta'}\subset\R^{2n-1}\subset Q$ with center $(c_1,0,\cdots,0)$ such that all boundary points $(\tilde x,f(\tilde x))$, $\tilde x\in V$, have a Levi form whose restriction to the corresponding complex tangent hyperplane is not positive semi-definite. Therefore, their complex tangent hyperplanes pass through the origin; the same for the real tangent hyperplane means that
$$
\sum_{j=1,j\neq k}^{2n}\frac{\partial f}{\partial x_j}(\tilde x)x_j=f(\tilde x),\quad \tilde x\in V.
$$
Using the Euler differential equation we see that
$$
f(t\tilde x)=tf(\tilde x),\quad \tilde x\in V,\; t\in I(\tilde x):=(1-\nu(\tilde x),1+\mu(\tilde x)),
$$
where $I(\tilde x)\tilde x\subset V$. In particular, we have for all boundary points $x=(\tilde x,f(\tilde x))$ sitting in $V\times I$ that $I(\tilde x)x\in\partial D$.

Take as a new defining function of $D$ near the point $c$ the following one: $\rho(z):=x_k-f(x_1,\dots,x_{k-1},x_{k+1},\dots,x_n,y_1,\dots,y_n)$.
Observe now, using homogeneity of $f$ near $(c_1,\tilde 0)$, that its gradient is constant near $c$ along the boundary of $D$. Call its norm by $\alpha$.
Then, again by homogeneity of $f$, we see that for $\tau\in (t_0-\eta,t_0]$, $0<\eta$ small, we have $(t_0-\eta,1]\subset\partial D$ and
$$
s(\tau b)=s_\rho(\tau b)/\alpha=\frac{t_0}{\alpha\tau} s_\rho(t_0b)=s(b)/\tau,
$$
where $s_\rho$ is defined as $s$ before but now with respect to the defining function $\rho$;
a contradiction to the definition of the infimum.

To summarize what we know so far (*): if  $b\in\partial_gD$, $b\neq 0$, and $\mathcal Lr(b;\cdot)$, restricted to the complex tangent hyperplane $T_b^\C(\partial D)$, is not positive semi-definite, then  $t_b<1$ and either $t_b=0$ or if $t_b>0$, then $t_bb\in\partial D\setminus\partial_gD$.

Now fix such a point, say again $b$; we may assume that $b=(b_1,\dots,b_n)$ with $b_1=\Re b_1<0$. From above we know that the segment $(t_bb,b]$ belongs to $T_b^\C(\partial D)$. Then the sphere $S_b$ with center $0$ passing through
$b$ intersect $\partial D$ transversally at $b.$ So, there exists a neighborhood $U$ of $b$ such that
$P_b=\partial D\cap S_b\cap U_b$ is a smooth surface of real codimension 2 containing only $\mathcal C^2$-smooth
non-pseudoconvex boundary points. Fix a $t\in(0,1)$ and put $Q_b^t=\{c\in P_b:t_c\ge t\}.$ By (*) the set $R_b^t=\{t_cc:c\in Q_t\}$
(possible empty) contains only non $\mathcal C^2$-smooth boundary points and hence $\mathcal H^{2n-2}(R_b^t)=0.$ Note that if
$(\B_n(x_k,r_k))_k$ is a covering of $R_b^t$ with $m:=\inf_k||x_k||>0,$ then $(\B_n(y_k,s_k))_k$ is a covering of $Q_b^t$
where $y_k=\frac{||c||}{||x_k||}x_k$ and $s_k=\frac{2}{m}r_k.$ Hence $\mathcal H^{2n-2}(Q_b^t)=0$. Therefore, there is a dense subset $P_b^t$ of $P_b$ such that
$(t,1]c\in\partial D$ for any $c\in P_b^t$. In particular, $(t,1]b\in\partial D$. And because of the arbitrariness of $t$ we end up with the information, that $(0,1]b\in\partial D$ and so $0\in\partial D$ and $b\in T_0^\C(\partial D)$.

Observe that $b$ was arbitrarily chosen. Therefore, every $\mathcal C^2$-boundary points $a\in\partial_g D$ at which $D$ is not pseudoconvex lie on the tangent plane at $0$ to $\partial D$. Fix one such $a_0$. Then this set (of $a$'s) is locally near $a_0$ an open part of a real hypersurface which is contained in $T_0^\R(\partial D)$. So $\partial D$  is locally near $a_0$ equal to this linear real hyperplane, which implies that $a_0$ is a $\mathcal C^2$-boundary point with vanishing Levi form on $T_{a_0}^\C(\partial D)$; a contradiction.
\end{proof}

\begin{remark}{\rm (a) Let $s$ be the function which was introduced during the proof of Theorem \ref{add}. If we replace the assumption of $\mathcal C^1$-smoothness by boundedness of $s$ from below, the result remains true and the proof above becomes essential easier.

(b) The proof above also shows that the $C^1$-smoothness can be replaced by the weaker condition that the convex hull of the
tangent cone (in sense of Whitney) at any boundary point is of real codimension at least 1.}
\end{remark}

Using an old result from \cite{Gra-Rem} we also have the following consequence.

\begin{corollary} Let $D\subset\C^n$ be a $\mathcal C^1$-smooth open set such that the set $A$ of non $\mathcal C^2$-smooth boundary points of $D$ has zero
$(2n-2)$-Hausdorff measure and it is locally contained in an analytic subset of codimension at least $1$ (any discrete set has these properties). Then $D$ has no
exceptional points.
\end{corollary}

\begin{remark} {\rm Finally, let us pose the following problems.}

\begin{enumerate}
\item Does an arbitrary non-pseudoconvex set $D\subset\C^3$ allow a complex line $l$ such that $\C^3\setminus S\subset l$ {\rm (see footnote 1)}?
\item Does the exceptional complex line in Theorem \ref{2} have more than one exceptional point?
\item Let $D\subset\C^n$ be a $\mathcal C^1$-smooth non-pseudoconvex open set such that $\mathcal H^{2n-2}(\partial D\setminus\partial_g D)=0$. Is there automatically a
non-pseudoconvex boundary point sitting in $\partial_g D$ {\rm (}compare with \cite{Zai-Zam}{\rm )}?
\end{enumerate}

\end{remark}

\end{document}